\documentclass{elsart}

\usepackage{amsmath,amsfonts,amssymb,amscd}

\newtheorem{defi}[thm]{Definition}
\newtheorem{rek}[thm]{Remark}

\newcommand{\twocase}[5]{#1 \begin{cases} #2 & \text{{\rm #3}}\\ #4
&\text{{\rm #5}} \end{cases}   }

\newcommand{\ncr}[2]{{#1 \choose #2}}
\newcommand{\gep}{\epsilon}

\newcommand\be{\begin{equation}}
\newcommand\ee{\end{equation}}
\newcommand\bea{\begin{eqnarray}}
\newcommand\eea{\end{eqnarray}}
\newcommand\bi{\begin{itemize}}
\newcommand\ei{\end{itemize}}
\newcommand\ben{\begin{enumerate}}
\newcommand\een{\end{enumerate}}
\newcommand\bc{\begin{center}}
\newcommand\ec{\end{center}}
\newcommand\ba{\begin{array}}
\newcommand\ea{\end{array}}

\newcommand{\Z}{\mathbb{Z}}

\newcommand{\ga}{\alpha}     

\newcommand{\jst}[1]{{#1\overwithdelims () p_2}}
\newcommand{\jsm}[1]{ { \underline{#1} \choose m} }            

\newcommand{\om}{\omega}
\newcommand{\h}[1]{\frac{#1}{2}}
\newcommand{\mc}[1]{\mathcal{#1}}
\newcommand{\symd}{\bigtriangleup}
\newcommand{\size}[1]{\left|#1\right|}

\begin{document}

\begin{frontmatter}




\title{Incomplete Quadratic Exponential Sums in Several Variables}

\author{Eduardo Due\~nez}\ead{eduenez@math.utsa.edu}\address{Department of Applied Mathematics,
The University of Texas at San Antonio, San Antonio, TX, 78249}
\author{Steven J. Miller}
\ead{sjmiller@math.brown.edu}
\address{Department of Mathematics, Brown University, Providence, RI 02912}
\author{Amitabha Roy}\ead{aroy@cs.bc.edu}\address{Department of Computer Science, Boston
College, Chestnut Hill, MA 02467}
\author{Howard Straubing}\ead{straubin@cs.bc.edu}\address{Department of Computer
Science, Boston College, Chestnut Hill, MA 02467}

\begin{abstract} We consider incomplete exponential sums in several
variables of the form
\begin{equation}  S(f,n,m) \ = \ \frac{1}{2^n} \sum_{x_1 \in \{-1,1\}} \cdots
\sum_{x_n \in \{-1,1\}} x_1 \cdots x_n\  e^{2\pi i
f(x)/p},\nonumber\
\end{equation}
where $m>1$ is odd and $f$ is a polynomial of degree $d$ with
coefficients in $\Z/m\Z$. We investigate the conjecture,
originating in a problem in computational complexity, that for
each fixed $d$ and $m$ the maximum norm of $S(f,n,m)$ converges
exponentially fast to $0$ as $n$ tends to infinity; we also
investigate the optimal bounds for these sums. Previous work has
verified the conjecture when $m=3$ and $d=2$. In the present paper
we develop three separate techniques for studying the problem in
the case of quadratic $f$, each of which establishes a different
special case.  We show that a bound of the required sort holds for
almost all quadratic polynomials, the conjecture holds for all
quadratic polynomials with $n \le 10$ variables (and the
conjectured bounds are sharp), and for arbitrarily many variables
the conjecture is true for a class of quadratic polynomials having
a special form.
\end{abstract}

\begin{keyword} incomplete exponential sums, boolean circuits

\MSC 11L07 (primary), 11G25 (secondary)
\end{keyword}
\end{frontmatter}


\section{Introduction}

We study sums of the form
\begin{equation}\label{eq:sfnm}  S(f,n,m) \ = \ \frac{1}{2^n} \sum_{x_1 \in \{-1,1\}} \cdots
\sum_{x_n \in \{-1,1\}} x_1 \cdots x_n\  \omega^{f(x)},
\end{equation}
where $m>1$ is odd, $\omega=e^{2\pi i/m}$, and $f$ is a polynomial
with coefficients in $\Z/m\Z$. This is an incomplete exponential
sum as each $x_i$ ranges only over $\{-1,1\}$.

Let $d$ be the degree of $f$.  It has been conjectured (see
\cite{green4,green3}) that there exists a positive $c_{m,d}<1$
such that \begin{equation}\label{eq:conjsfnm} |S(f,n,m)|\ \le \
c_{m,d}^n.\end{equation} Exponential sums have a rich history, and
estimates of their size have numerous applications, ranging from
uniform distribution to solutions to Diophantine equations to
$L$-functions to the Circle Method, to name a few. Our problem
originates in computer science, where \eqref{eq:sfnm} arises in
the study of the complexity of boolean circuits. The conjecture
\eqref{eq:conjsfnm} implies that a very special kind of $n$-input
boolean circuit, containing ``mod-$m$ gates''---that is, gates
that determine whether the number of their input bits that are on
is divisible by $m$---requires exponentially many (in $n$) gates
in order to simulate a single mod-2 gate (i.e., in order to
``compute parity''). Such questions concerning exponential lower
bounds on the size of circuits that perform various computations,
and, in particular, the relation between the computing power of
modular gates with different moduli,  are notoriously difficult,
and progress in this area has been quite scant.  See Green
\cite{green} for a precise account of the connection between this
problem and circuit complexity.

It is known (Alon and Beigel \cite{alonbeigel}) that for each
fixed $n$, $d$ and $m$ there exists a positive constant
$b_{d,m,n}$ such that
\begin{equation}|S(f,n,m)|\ < \ b_{d,m,n},\end{equation}and
\begin{equation}\lim_{n\to\infty}b_{d,m,n}=0.\end{equation} This theorem is
proved using Ramsey-theoretic techniques, and the resulting
sequences converge very slowly to 0. In terms of computational
complexity, this only tells us that the minimum circuit size
required to compute parity of $n$ bits tends to infinity with $n$.
It is of far more interest, from the computational point of view,
to show exponentially fast growth in minimum circuit size.  This
is generally interpreted as showing that parity circuits of the
required kind cannot feasibly be built.

The conjecture \eqref{eq:conjsfnm} holds trivially for $d=1$,
since in this case $S(f,n,m)$ is a product of a complex number of
norm 1 and $n$ factors of the form $\omega^k-\omega^{-k}$.  In the
case $d=2$, \eqref{eq:conjsfnm} has been proved only in the case
$m=3$, and the optimal value of $c_{3,2}$ determined (see
 \cite{green}); however this proof appears to shed no light on
what occurs with other odd moduli. The conjecture has also been
verified (see \cite{green3}) when $f$ is a symmetric polynomial in
$n$ variables, of  poly-logarithmic degree (in $n$) and  for any
odd modulus $m$.

A natural approach to proving \eqref{eq:conjsfnm} is to use
Weil-type bounds for multiple exponential sums. While there have
been many bounds published for incomplete and complete exponential
sums over many variables (see Notes to Chapter 5 of \cite{lidl},
as well as
\cite{chubarikov,chub,davlew,deligne,loxton,mord,tiet}), none
seems to apply to our situation so far. We quickly review these
approaches; the inapplicability of these techniques led us to the
methods of this paper.

Consider the bounds of incomplete exponential sums from
\cite{mord,tiet} with $m$ an odd prime $p$. Though not directly
applicable to our problem because of the factor $x_1\cdots x_n$,
it is enlightening to see what bounds estimates of this type can
generate. Using finite Fourier transforms, these represent the
incomplete sum as $\frac{2^n}{p^n}$ times a complete sum plus an
error term. The bounds for the error term are improved if we are
summing over consecutive $x_i$ (this can readily be done for our
problem by sending $x_i$ to $\frac{x_i+1}2$; the factor $x_1\cdots
x_n$ is replaced with $2^n$ terms, but each term is divided by an
additional factor of $2^n$). For example, Mordell \cite{mord}
considers incomplete sums \begin{equation}S_n' \ = \ \sum_{0\le
x_1 < \ell_1} \cdots \sum_{0 \le x_n < \ell_n} e_p(f(x)), \ \ \
e_p(x) \ = \ e^{2\pi i x/p}. \end{equation}Denote the complete sum
by $S_n$. If $t = (t_1,\dots,t_n)$ has $r$ non-zero entries,
suppose there is a constant $E_n^{(r)}$ (independent of $t$ but
depending on $p$ and $f$) such that
\begin{equation}\left|\sum_{x_1 \bmod p} \cdots \sum_{x_n \bmod p}
e_p(f(x) + t_1x_1 + \cdots + t_n x_n)\right| \ \le \ E_n^{(r)};
\end{equation}In general we expect $E_n^{(r)}$ to be at least $p^{n/2}$.
Mordell proves that \begin{equation}S_n' \ = \ \frac{\ell_1\cdots
\ell_n}{p^n} S_n + \Theta_n^{(n)} E_n^{(n)} \log^n p + R_n,
\end{equation}where $|\Theta_n^{(n)}| < 1$ and
\begin{equation}R_n \ = \ \sum_{r=1}^{n-1} \frac{\ell_{r+1}\cdots
\ell_n}{p^{r-n}}\ \Theta_n^{(r)} E_n^{(r)} \log^r p, \ \ \
|\Theta_n^{(r)}| < 1. \end{equation}For $p > 3$, the bounds for
$E_n^{(r)}$ are too weak. The reason for the failure of these
methods is the paucity of points in the sub-variety we sum over;
we would need to let the number of $x_i$ we sum over grow with
$p$.

It is possible to transform our incomplete exponential sum to a
complete one involving Legendre symbols by having the variables
range over all of $\Z/m\Z$ (this was already observed by
\cite{green}, however we show an alternate method here). For ease
of exposition we assume now that $m$ is an odd prime congruent to
$-1$ modulo $4$. In this case, $\jsm{-1} = (-1)^{(m-1)/2} = -1$
and we have
\begin{equation}\begin{split}
 S(f,n,m) & \ = \
\frac{1}{2^n} \sum_{x_1 \in \{-1,1\}} \cdots \sum_{x_n \in
\{-1,1\}}
x_1 \cdots x_n e_m\left( f(x) \right)\\
& \ = \ \frac{1}{2^n} \sum_{x_1 \in \{-1,1\}} \cdots \sum_{x_n \in
\{-1,1\}} x_1^{(m-1)/2} \cdots x_n^{(m-1)/2}\\ & \ \ \ \ \ \ \ \ \
\ \ \ \ \ \ \times \ e_m\left(
f(x_1^{(m-1)/2},\dots,x_n^{(m-1)/2}) \right).\end{split}
\end{equation}
The above weakly depends on $x_i$; all that matters is the value
of $\jsm{x_i}$, the Legendre symbol. Thus we may extend all
summations from $x_i \in \{-1,1\}$ to $x_i \in \Z/m\Z$ (note we
may trivially include any $x_i = 0$). Letting $g(x) =
f(x_1^{(m-1)/2},\dots,x_n^{(m-1)/2})$ we are led to a new
formulation of the problem. Namely, we must estimate
\begin{equation}\label{eq:Sgnm} S(g,n,m) \ = \ \frac{1}{(m-1)^n} \sum_{x_1 =
0}^{m-1} \cdots \sum_{x_n = 0}^{m-1} \jsm{x_1} \cdots \jsm{x_n}
e_m\left( g(x) \right). \end{equation}This is a mixed exponential
sum, involving multiplicative (the Legendre symbol) and additive
(the exponential function) characters. When there are no Legendre
symbols in \eqref{eq:Sgnm}, one often obtains bounds of the form
\begin{equation}\label{eq:boundsfordmn} (d-1)^n m^{n/2},\end{equation}where $d$
is the degree of the highest homogeneous component, $m$ is the
modulus, and $n$ the number of variables (see \cite{deligne}). The
substitution (replacing $f$ with $g$) increases the degree $d$ too
much for the general Weil-Deligne type bounds to help, except when
$m=3$ where the conjecture is already known. Note the degree of
$g$ is $m-1$, so the degree increases unless $m=3$. For $m=3$ this
does lead to a new proof of the conjecture for special $f$ (see
Appendix \ref{sec:appm3d2} for details).

An alternate approach to \eqref{eq:sfnm} is to rewrite it as
\begin{eqnarray}\label{eqvarietysum}& &  \frac{1}{m^n} \frac{1}{2^n}
\sum_{\ga_1,\dots,\ga_n \bmod m} \sum_{x_1,\dots,x_n \in \{-1,1\}
}x_1\cdots x_n \left[ \prod_i e_m( \ga_i (x_i^2 - 1) ) \right]
\nonumber\\ & & \ \ \ \ \  \ \ \ \ \  \ \ \ \ \  \ \ \ \ \  \ \ \
\ \  \ \ \ \ \  \ \ \ \ \  \ \ \ \ \ \times \ e_m(
f(x_1,\dots,x_n) ). \end{eqnarray} In the bracketed product, the
sum over each $\ga_i$ is $0$ unless $x_i^2 - 1 \equiv 0$ mod $m$;
in other words, we may extend the summation over each $x_i$ to be
over all of $\Z/m\Z$. Note it is relatively easy to explicitly
incorporate summing over the sub-variety $x_i^2 = 1$.
Unfortunately, the number of variables of the new polynomial is
now $2n$, and the degree is now $3$. This will also be a poor
substitution. Again ignoring the $x_1\cdots x_n$, the bounds from
\eqref{eq:boundsfordmn} are of the form
\begin{equation}\frac{1}{m^n \cdot 2^n} \cdot (3-1)^{2n} m^{2n/2}
\ = \ 2^n, \end{equation}which is too large; other similar bounds
also just fail (see for example \cite{chubarikov}).

In the present paper we investigate the sums $S(f,n,m)$ from
\eqref{eq:sfnm} in the case $d=2$ and arbitrary odd $m$.  In this
setting the conjecture takes on a sharper form, since we believe
we know the optimal value of $c_{m,2}$ and the quadratic
polynomials $f$ for which the optimal bound is attained. While we
have not settled the question, we have developed three quite
different techniques for studying the problem. Each of these
methods produces a proof of a different special case of the
conjecture for quadratic polynomials. We believe that at least one
of these methods, or some combination of them,  can be pushed
further to settle the general problem.

We first investigate the conjecture probabilistically by
evaluating the higher-order moments of $|S(f,n,m)|$ as $f$ ranges
over the set of all quadratic polynomials in $n$ variables.  As a
result, we are able to show that if $\gamma<1$ is quite close to
1, then all but an exponentially small (in $n$) proportion of the
$|S(f,n,m)|$ are bounded by $\gamma^n$.

We then give a detailed analysis of the structure of these sums
for small $n$. As a consequence, we are able to prove our
conjectured upper bound holds whenever $n\leq 10$ for any odd $m$.
Further, we prove these bounds are sharp for $n\le 10$.

Finally, we interpret $S(f,n,m)$ as a coefficient in the Fourier
expansion of $\omega^{f(x_1,\ldots,x_n)}$, when this function is
viewed as an element of $L^2(\{-1,1\}^n)$. We are able, for a
large class of polynomials, to determine the Fourier expansion
directly, and thus obtain the conjectured bound.

\section{Definitions and Statement of Main Results}\label{sec:defnot}
\setcounter{equation}{0}

Let $m$ be a fixed odd integer and let $f(x) = f(x_1,\dots,x_n)
\in \Z[x_1, x_2, \ldots, x_n]$ be a polynomial with integral
coefficients of degree at most 2 in $n$ variables. We are
interested in finding sharp upper bounds to the norm of
\begin{equation}  S(f,n,m) \ = \ \frac{1}{2^n} \sum_{x_1 \in \{-1,1\}} \cdots
\sum_{x_n \in \{-1,1\}} x_1 \cdots x_n\  \om^{f(x)},
\label{defeqn}
\end{equation}
where $\om \ = \ e^{2\pi i / m}$ is the principal $m$-th root of
unity. Letting $e_m(z) = e^{2\pi i z/m}$, we often write
$\om^{f(x)} = e_m(f(x))$. When $n$ and $m$ are obvious from the
context, we refer to this sum as $S(f)$.  These are incomplete
exponential sums, as each $x_i$ is restricted to lying in
$\{-1,1\}$; the easier case has each $x_i \in \Z/m\Z$. It is
important to note that for our applications, the modulus $m$ is
fixed and our goal is to study the norm of the $S(f, n, m)$ as $n$
and $f$ vary. We shall refer to $S(f,n,m)$ as the
\emph{normalized} sum, on occasion referring to the unnormalized
sum $2^n S(f,n,m)$ as $\tilde{S}(f,n,m)$. The philosophy of
square-root cancellation suggests that $\tilde{S}(f,n,m)$ should
typically be of size $2^{n/2}$.

Without loss of generality, we may assume there are no diagonal or
constant terms in $f(x)$: as each $x_i \in \{-1,1\}$, $x_i^2$ is
constant and hence does not affect $\size{S(f)}$. Thus we restrict
our attention to  $f(x)$ of the form \begin{equation}f(x) \ = \
\sum_{1 \le i < j \le n} a_{ij} x_i x_j \ + \ \sum_{1 \le k \le n}
b_k x_k. \end{equation}and we refer to this set of polynomials as
$\Z_m^2[x_1, x_2, \ldots, x_n]$, or $\Z_m^2[n]$ for short.

For fixed $n$ and $m$, let $\mathcal{F} \subset \Z_m^2[n]$ be an
arbitrary family of polynomials. For $r
>0$, we define the $r$\textsuperscript{th} moment of $\mathcal{F}$,
 denoted by $M_{r,\mathcal{F}}$, by
 \begin{equation}
 M_{r, \mc{F}}  \ = \  \langle |S(f,n,m)|^{r} \rangle_\mathcal{F}
 \ = \  \frac{1}{\size{\mc{F}}}\ \sum_{f \in \mc{F}}
|S(f,n,m)|^{r}.
\end{equation}
When $\mc{F}$ is obvious from the context, we write $M_r$ for the
$r$\textsuperscript{th} moment.

We now define a few parameters that appear in our results:
\begin{itemize}
\item\noindent $c:=\lfloor{{m+1}\over 4}\rfloor \in \Z$. This
value  maximizes $|\omega^y-\omega^{-y}|$. \item  \noindent $q:=
\size{\omega^c - \omega^{-c}}=2\cos{{\pi}\over{2m}}$. \item
\noindent $r := \cos{{3\pi}\over{2m}}$ \   denotes the second
largest value of $|\omega^y-\omega^{-y}|$.  A simple calculation
shows that this is attained when $y=\lfloor{{m+3}\over 4}\rfloor$.
\item  \noindent $s:=\cos{{\pi}\over{m}}$. This is the second
largest value of $|\omega^y+\omega^{-y}|$ (the largest value is 2,
when $y=0$).
\end{itemize}

Associated with every polynomial $f=\sum_{i,j} a_{ij} x_i x_j +
\sum_i b_i x_i \in \Z_m^2[n]$ (of degree $ \leq 2$) is an
undirected graph $G=G(f)$ with vertex set $\{1, \ldots, n\}$ and
edge set $\{\{i,j\}:a_{ij}\neq 0\}$. Recall that a tree is a
connected acyclic graph and a forest is a collection of trees.

Our main result towards proving the conjectured bounds in
\eqref{eq:conjsfnm} is

\begin{thm}\label{grandmother}\
\begin{itemize}
\item[(i)] Let $\mc{F}$ (resp. $\mc{G}$) denote the set of all
quadratic polynomials (resp. homogeneous quadratic polynomials) in
$\Z_m^2[n]$. Then the second moments are \begin{equation}M_{2,
\mc{F}}\ =\ \frac{1}{2^n}, \ \ \ \ \ M_{2, \mc{G}} = \
\frac{1+(-1)^n}{2^n}. \end{equation}Furthermore, for $m > 3$, the
sixth moment satisfies \begin{equation}M_{6, \mc{F}} \ \le \
\frac{9n(n-1) + (9n+1) 2^{2-2n} }{4} \frac{1}{2^{3n}}.
\end{equation}\item[(ii)] For all odd $m \geq 3$ and $n \leq 10$,
\be|S(f,n,m)|\ \leq \ \left({q\over 2}\right)^{\lfloor{{n+1}\over
2}\rfloor}.\end{equation}This bound is sharp, as there are
polynomials where equality holds. \item[(iii)] If $f \in
\Z_m^2[n]$ is such that $G(f)$ becomes  a forest of trees on
deletion of at most $(n-2) \log(2/q)$  edges from $G(f)$, then
 \begin{equation}\size{S(f,n,m)}\ \leq\  \biggl({q\over 2}\biggr)^{\lfloor{{n+1}\over 2}\rfloor}. \end{equation}Additionally, if $G(f)$ is itself a  tree, then
 \begin{equation}\size{S(f,n,m)}\ \leq\ \left(q \over 2 \right)^{n-1} .\end{equation}\end{itemize}
\end{thm}

The moment bounds in Theorem \ref{grandmother}~(i) allows us to
estimate the number of polynomials with large norms. Specifically,
we prove:

\begin{cor}\label{child}
 Let $f \in \Z_m^2[n]$ be chosen randomly and uniformly from $\Z_m^2[n]$. Then
for any $\gamma > 0$, \begin{equation}\frac{\frac{1}{2^n} -
\gamma^{2n}}{1 - \gamma^{2n}} \ \le \ \text{\rm
Prob}\left(|S(f,n,m)|\ \ge\ \gamma^n\right) \ \le \
\min\left(\frac{1}{(2\gamma^2)^n},
\frac{9n(n+1)/4}{(2\gamma^2)^{3n}}\right).
\end{equation}\end{cor}

\begin{rek}  A critical case occurs when $\gamma = \frac{1}{\sqrt{2}}$.
This occurs when we have square-root cancellation. The second and
sixth moment bounds, at $\gamma = \frac{1}{\sqrt{2}}$, give no
information: $0 \le P(\gep) \le 1$. In other words, we cannot
obtain more than square-root cancellation on a positive proportion
of polynomials. This agrees nicely with the philosophy that
square-root cancellation is the best one can hope for in general.
\end{rek}

The previous remark yields the following negative result:

\begin{cor} For any $\gamma < \frac{1}{\sqrt{2}}$, at least an
exponentially small (in $n$) proportion of the $f$, independent of
$m$, satisfy $|S(f,n,m)| \ge \gamma^n$.
\end{cor}

The bounds in Theorem \ref{grandmother} and ample experimental
evidence for small values of $n$ lead us to make the following
conjecture:

\begin{conj}\label{holygrail}
Let $m \geq 3$ be odd and let $n$ be a non-negative integer. Then
\begin{equation}\size{S(f,n,m)}\ \leq\ \left( q \over 2 \right)^{\lfloor
\frac{n+1}{2} \rfloor }.\end{equation}Moreover, the upper bound is
attained by all polynomials of the form \begin{equation} c (\pm
x_1 x_2 \pm x_3 x_4\pm \cdots \pm x_{n-1} x_n) \end{equation}when
$n$ is even, and by any polynomial of the form \begin{equation}c
(\pm x_1 x_2 \pm x_3 x_4\pm \cdots \pm x_{n-1} x_n \pm x_{n+1})
\end{equation}when $n$ is odd, where the constant $c = \lfloor (m+1)/4
\rfloor$.
\end{conj}

Note that the special case of Conjecture \ref{holygrail} has
already been verified for all $n$ and $m=3$ \cite{green}. Green's
proof for $m=3$ makes use of special relations that hold between
the third roots of unity, and we have not been able to generalize
these equations to higher roots. \\

\noindent\emph{Organization of paper}: We prove Theorem
\ref{grandmother}(i) in Section \ref{momentsection}, Theorem
\ref{grandmother}(ii) in Section \ref{smallvarsection} and finally
in Section \ref{fouriersection} we prove Theorem
\ref{grandmother}(iii). In Section \ref{sec:addapproachfuturework}
we discuss a generalization of Conjecture \ref{holygrail} and
future work.

\section{Bounds through Moments} \label{momentsection}

In this section, we prove Theorem \ref{grandmother}(i) and
Corollary \ref{child} by computing the moments of the exponential
sums $S(f,n,m)$. We can compute the second moment exactly, while
for the sixth moment we provide an upper bound. These calculations
enable us to provide estimates on the proportion of polynomials
with large norm. Theorem \ref{grandmother}(i) follows immediately
from Theorems \ref{thm:2ndnonhomog}, \ref{thm:2ndhomog} and
\ref{thm:6thnonhomog}, while Corollary \ref{child} follows from
Theorem \ref{grandmother}(i) and Theorem
\ref{thm:boundsfrommoments}.

\subsection{Moment Bounds}

Using moments, one can gain information about the maximum value of
$|S(f,n,m)|$. As $r\to\infty$, the $r$\textsuperscript{th} root of
the $r$\textsuperscript{th} moment converges to the largest value
of $|S(f,n,m)|$. Unfortunately, because of combinatorial
complications, we cannot compute high enough (in $n$) moments to
obtain the desired bounds for individual $S(f,n,m)$, as the order
of the moment needed tends to infinity with $n$. Thus, while the
method of moments allows us to conclude that ``most'' $S(f,n,m)$
have the desired cancellation, to obtain these estimates for all
$S(f,n,m)$ requires, at present, moments that are too
combinatorially difficult to calculate. We do observe that the low
moments are growing at a rate which is indicative of the
conjectured bounds being true.

\begin{defi}[$P(\gep)$]
\begin{equation}P(\gep) \ = \ \mbox{\rm Prob}\left(|S(f,n,m)| \ge \gep\right).
\end{equation}\end{defi}

\begin{thm}[Bounds from Moments]\label{thm:boundsfrommoments}
Assume $L_r \le M_r \le U_r$. Then \begin{equation} \frac{L_r -
\gep^r}{1 - \gep^r} \ \le \ P(\gep) \ \le \  \frac{U_r}{\gep^r}.
\end{equation}
\end{thm}

\textbf{Proof.} As
 \begin{equation}0^r \cdot (1 - P(\gep)) + \gep^r \cdot P(\gep) \ \le \ U_r,
\end{equation}we obtain  \begin{equation}P(\gep) \ \le \ \frac{U_r}{\gep^r}.
\end{equation}The above is just Chebychev's Inequality, which allows us to
measure the ``bad'' set of $f$. The lower bound follows from
\begin{equation}\gep^r \cdot \left(1 - P(\gep) \right) + 1 \cdot
P(\gep) \ \ge \ L_r.  \end{equation} \qed

Good bounds can be found for any fixed moment (if one is willing
to do enough work); we provide details for the second moment
(which is very straightforward) and the sixth moment (which
illustrates the type of complications that arise in studying the
higher moments).

We now bound the second and sixth moments. Recall $e_m(x) =
e^{2\pi i x / m}$. We constantly use the following observation:

\begin{lem}\label{lem:expsumgs} For any positive integer $m$, \begin{equation}\twocase{\sum_{a \bmod m}
e_m(ar) \ = \ }{m}{if $r \equiv 0 \bmod m$}{0}{otherwise.}
\end{equation}\end{lem}

\textbf{Proof.} If $r\equiv 0 \bmod m$, each term is 1 and the
claim is clear. Otherwise the above is a geometric series with
ratio $e_m(r)$, equal to $\frac{e_m(0r) - e_m(mr)}{1-e_m(r)}=0$.
\qed

\subsection{The Second Moment}

\subsubsection{All Quadratic Polynomials in $\Z_m^2[n]$}

\begin{thm}\label{thm:2ndnonhomog} Let $\mathcal{F} = \Z_m^2[n]$.
Then for any integer $m \ge 2$, \begin{equation}M_2 \ = \
\frac1{2^n}. \end{equation}\end{thm}

\textbf{Proof.} The second moment of $|S(f,n,m)|$  is
\begin{equation}\begin{split} M_2
 & \ = \  \frac{1}{|\mathcal{F}|} \sum_{a_{ij} \bmod m}
\sum_{b_k \bmod m} \left(\frac{1}{2^n} \sum_{x_1 \in \{-1,1\}}
\cdots \sum_{x_n \in  \{-1,1\}} x_1 \cdots x_n e_m\left( f(x)
\right) \right) \\ & \ \ \ \ \ \ \ \ \ \ \ \ \ \ \cdot
\left(\frac{1}{2^n} \sum_{y_1 \in \{-1,1\}} \cdots \sum_{y_n \in
\{-1,1\}} y_1 \cdots y_n e_m\left( -f(y) \right)
\right).\end{split}\end{equation} Interchanging summations, for a
fixed $2n$-tuple $(x_1,\dots,y_n)$, we have terms such as
\begin{equation}\sum_{a_{ij} \bmod m} \sum_{b_k \bmod m} e_m\left(
f(x) - f(y) \right). \end{equation}This equals
\begin{equation}\sum_{a_{ij} \bmod m} \sum_{b_k \bmod m} e_m\left(
\sum_{i,j} a_{ij} (x_i x_j - y_i y_j) + \sum_k b_k (x_k - y_k)
\right). \end{equation}If $x_k \not\equiv y_k \bmod m$, then by
Lemma \ref{lem:expsumgs} the sum over that $b_k$ is zero. Thus the
only non-zero contributions for a $2n$-tuple are when each $x_k$
equals the corresponding $y_k$. There are $2^n$ such tuples. Note
that in this case, each sum over $b_k$ gives $m$. Further, each
sum over an $a_{ij}$ also gives $m$, as $x_i x_j - y_i y_j \equiv
0 \bmod m$.

Thus for each of the $2^n$ tuples where $x_k = y_k$, the sums over
$a_{ij}$ and $b_k$ give $m^{n(n+1)/2} = |\mathcal{F}|$, and
$x_ky_k = 1$. Substituting yields \begin{equation} M_2 \ = \
\frac{1}{|\mathcal{F}|} \cdot \frac{1}{2^{2n}} \cdot 2^n \cdot
|\mathcal{F}| \ = \ \frac{1}{2^n}.
\end{equation}
\qed

\begin{rek} Theorem \ref{thm:2ndnonhomog}
implies that  on average there is square-root cancellation; using
the Cauchy-Schwartz inequality, we find \begin{equation}\langle
|S(f,n,m)|\rangle_{\mathcal{F}} \ \le \
\frac1{2^{n/2}}.\ee\end{rek}

\subsubsection{Homogeneous Quadratic Polynomials in $\Z_m^2[n]$}

While we are primarily interested in bounds for $S(f,n,m)$ for
non-homogeneous $f$, we quickly investigate the homogeneous case.

\begin{thm}\label{thm:2ndhomog} Let $\mathcal{G}$ be the family of all
homogeneous quadratic polynomials in $\Z_m^2[n]$. Then
\begin{equation}M_2 \ = \ \frac{1+(-1)^n}{2^n}. \end{equation}\end{thm}

\textbf{Proof.} As this case is similar to the previous one, we
just sketch the arguments below. The main difference is we now
only have sums over $a_{ij} \bmod m$; there are no $b_k$ sums.
Thus for each $2n$-tuple $(x_1,\dots, y_n)$, we have factors such
as \begin{equation}\sum_{a_{ij} = 0}^{m-1} e_m\left( a_{ij}(x_i
x_j - y_i y_j) \right). \end{equation}If $x_i x_j - y_i y_j \equiv
0 \bmod m$ then the $a_{ij}$-sum is $m$; otherwise, it is $0$. As
$m$ is odd, if $x_ix_j - y_iy_j \equiv 0 \bmod 1$, then it equals
zero.

There are two possibilities. First, each $y_i$ could equal $x_i$.
Then clearly all relevant terms equal 0. For the second
possibility, assume there exists an $i$ such that $x_i = -y_i$.
Then for any $j \neq i$, $x_i x_j - y_i y_j = 0$ becomes $x_j +
y_j = 0$. Therefore, if one $y_i = -x_i$, then \emph{all} $y_i = -
x_i$. We again find the $a_{ij}$-sum equals $m$.

Therefore, for each $n$-tuple $(x_1,\dots,x_n)$ there are two
$y$-tuples, $(x_1,\dots,x_n)$ and $(-x_1,\dots,-x_n)$. The
exponential sums over $a_{ij}$ give $m^{n(n-1)/2} =
|\mathcal{G}|$. We then multiply by \begin{equation}x_1 \cdots x_n
x_1 \cdots x_n + x_1 \cdots x_n (-x_1) \cdots (-x_n)  = 1 +
(-1)^n, \end{equation}and find that \begin{equation} M_2 \ = \
\frac{1}{|\mathcal{G}|} \cdot \frac{1}{2^{2n}} \cdot 2^n \cdot
\left(1 + (-1)^n\right) \cdot |\mathcal{G}| \ = \ \frac{1 +
(-1)^n}{2^n}. \end{equation} \qed

Note if $n$ is odd, the second moment is 0, which implies that
$S(f,n,m) = 0$ for all $f$; this is also seen by comparing the
contributions from $(x_1,\dots,x_n)$ and $(-x_1,\dots, -x_n)$.

\subsection{The Sixth Moment}

\begin{thm}\label{thm:6thnonhomog} Assume $m > 3$ is odd. The sixth moment for
$\mathcal{F} = \Z_m^2[n]$ satisfies \begin{equation}M_6 \ \le \
\frac{9n(n-1) + (9n+1)2^{2-2n}}{4} \frac{1}{2^{3n}} \ \sim \
\frac{9n(n-1)}{4} \frac{1}{2^{3n}}. \end{equation}\end{thm}

\textbf{Proof.} We have six tuples in the calculation of the sixth
moment, say $X_1 = (x_{1,1},\dots,x_{1,n})$ to $X_6 =
(x_{6,1},\dots,x_{6,n})$. We have exponential factors such as
\begin{equation}\label{eq:aijexp} \sum_{a_{ij} \bmod m} e_m\Big(
a_{ij} (x_{1,i}x_{1,j} + x_{2,i}x_{2,j} + x_{3,i}x_{3,j} -
x_{4,i}x_{4,j} - x_{5,i}x_{5,j} - x_{6,i}x_{6,j}) \Big)
\end{equation}and
\begin{equation}\label{eq:bkexp}\sum_{b_k \bmod m} e_m\Big( b_k( x_{1,k} +
x_{2,k} + x_{3,k} - x_{4,k} - x_{5,k} - x_{6,k})\Big).
\end{equation}The $b_k$-sum is zero unless \begin{equation}\label{eqsixxs}
x_{1,k} + x_{2,k} + x_{3,k} - x_{4,k} - x_{5,k} - x_{6,k} \ \equiv
\ 0 \ \bmod \ m. \end{equation}
\begin{rek} If we were calculating the $2r$\textsuperscript{th} moment, we
would have \begin{equation}x_{1,k} + \cdots + x_{r,k} - x_{r+1,k}
- \cdots - x_{2r,k} \ \equiv \ 0 \ \bmod \ m. \end{equation}We
want to conclude that $x_{1,k} + \cdots - x_{2r,k} = 0$. As each
term is congruent to $1$ mod $2$, the sum is always even. For the
sixth moment, if the sum is congruent to zero mod $m$ then it is
zero unless $m = 3$; this is clear for $m > 6$, and if $m=5$ this
follows immediately. Thus some modifications are needed to use
these techniques for $m=3$; as the main theorem can be proved for
all $n$ for $m=3$, we do not explore such extensions here and
content ourselves with remarking that slight changes are needed
for small $m$ and larger moments (for example, $m=5$ and $2r =
12$). In all arguments below, we may replace congruent to $0$ mod
$m$ with equals $0$.
\end{rek}

Thus, in \eqref{eqsixxs}, if exactly $m$ of the first three
$x_{h,k}$'s are $+1$, then exactly $m$ of the last three
$x_{h,k}$'s are $+1$. For each $k$, there are four structurally
different ways to choose the $x_{h,k}$'s:

\ben

\item None of the $x_{1,k}, x_{2,k}, x_{3,k}$ are $1$; there is
$\ncr{3}{0}\ncr{3}{0} = 1$ way to do this.

\item Exactly one of the $x_{1,k}, x_{2,k}, x_{3,k}$ are $1$;
there are $\ncr{3}{1}\ncr{3}{1} = 9$ ways to do this.

\item Exactly two of the $x_{1,k}, x_{2,k}, x_{3,k}$ are $1$;
there are $\ncr{3}{2}\ncr{3}{2} = 9$ ways to do this.

\item Exactly three of the $x_{1,k}, x_{2,k}, x_{3,k}$ are $1$;
there is $\ncr{3}{3}\ncr{3}{3} = 1$ way to do this.

\een

We call these conditions (1) through (4). For all $(i,j)$, we have
\begin{equation}\label{eqsixxij} x_{1,i}x_{1,j} + x_{2,i}x_{2,j} +
x_{3,i}x_{3,j} - x_{4,i}x_{4,j} - x_{5,i}x_{5,j} - x_{6,i}x_{6,j}
\ = \ 0, \end{equation}or else the $a_{ij}$-sum is zero. We now
analyze the consequences of having one of the above conditions
hold.

For example, assume there is a $k_0$ such that condition~(1) holds
(all six of the $x_{h,k_0}$ are $-1$). Then for all $j \neq k_0$,
substituting into \eqref{eqsixxij} and multiplying through by $-1$
yields \begin{equation}x_{1,j} + x_{2,j} + x_{3,j} - x_{4,j} -
x_{5,j} - x_{6,j} \ = \ 0. \end{equation}This is exactly the
condition from the $b_k$-sums (\eqref{eq:bkexp} and
\eqref{eqsixxs}), and provides \emph{no} new information (ie, this
equation is already satisfied for all $j$). Thus, whenever
condition (1) is satisfied, no new information is obtained. In
effect, whenever condition (1) holds, it is as if we have a
smaller degree for our polynomial. This is primarily because
initially there are $2^6$ possibilities for a 6-tuple, and when
condition (1) holds, there is only one possibility.

Assume now condition~(2) holds for some fixed index $k_0$, namely
exactly one of the first three is $+1$, exactly one of the last
three is $+1$. There are 9 different ways this can occur; by
symmetry we can relabel so that $x_{1,k_0} = x_{4,k_0} = 1$.
Substituting into \eqref{eqsixxij} yields, for any $j \neq k_0$,
\begin{equation}\label{eq:sixnegj} x_{1,j} - x_{2,j} - x_{3,j} - x_{4,j} +
x_{5,j} + x_{6,j} \ = \ 0. \end{equation}However, from the
$b_k$-sum with $k = j$ (\eqref{eq:bkexp} and \eqref{eqsixxs}), we
have \begin{equation}\label{eq:sixposj} x_{1,j} + x_{2,j} +
x_{3,j} - x_{4,j} - x_{5,j} - x_{6,j} \ = \ 0.
\end{equation}Adding
 \eqref{eq:sixnegj} and \eqref{eq:sixposj} and dividing
by 2 (\emph{note here we use $m$ is odd!}) yields
\begin{equation}\label{eq:sixtwojone} x_{1,j} \ = \ x_{4,j},
\end{equation}while subtracting the two and dividing by 2 yields
\begin{equation}\label{eqsixtwoj} x_{2,j} + x_{3,j} \ = \ x_{5,j} + x_{6,j}.
\end{equation}
There are two possibilities in \eqref{eqsixtwoj}: we could have
each side is two equally signed summands, or oppositely signed
summands. We have already determined $x_{1,j} = x_{4,j}$; we now
isolate the relations among the other $x$'s in this case.

\begin{lem}\label{lem:sixj} Assume condition (2) holds for some
$k_0$, and for definiteness assume $x_{1,k_0} = x_{4,k_0}$. Then
for all $j\neq k_0$ we have $x_{1,j} = x_{4,j}$, and exactly one
of the following must hold:

\bi

\item If $x_{2,j} = x_{3,j}$, then $x_{2,j} = x_{3,j} = x_{5,j} =
x_{6,j}$. There are two ways this can occur (once the sign of
$x_{2,j}$ is chosen, all other values are determined). We call
this case ``equally signed terms''.

\item If $x_{2,j} = -x_{3,j}$, then $x_{5,j} = -x_{6,j}$. The two
possibilities are

\noindent \ \ \ $(i)$\ \ $x_{2,j} = -x_{3,j} = x_{5,j} =
-x_{6,j}$;

\noindent \ \ \ $(ii)$ $x_{2,j} = -x_{3,j} = -x_{5,j} = x_{6,j}$.

There are two ways for each possibility to occur; again, once
$x_{2,j}$ is chosen, the rest are determined. We denote this case
``oppositely signed terms''. \ei
\end{lem}

Note in all of the relations above, we always have $x_{1,j}\cdots
x_{6,j} = +1$; thus, the contributions from these terms will not
negatively reinforce. If there is some $k_0$ so that condition (2)
holds, then for each $j \neq k_0$, there are 12 choices for the
variables $(x_{1,j},\dots,x_{6,j})$, and each choice leads to a
contribution of $|\mathcal{F}|$. The reason there are 12 choices
is that there are two ways to satisfy $x_{1,j} = x_{4,j}$, and
then 6 ways to satisfy the other relations. There are $n$ ways to
choose an index $k_0$ such that condition (2) holds, and 9 ways to
choose the indices for that $k_0$. As there are $2^{6n} = 64^n$
6-tuples, this leads to condition (2) terms contributing at most
\begin{equation}9n \cdot \frac{12^{n-1}}{64^n} \ = \ \frac{9n}{12}
\left(\frac{12}{64}\right)^{n} \ = \ \frac{3n}{4}
\frac{1}{1.74716^{3n}}. \end{equation}For square-root
cancellation, the sixth moment should be of size $\frac1{2^{3n}}$;
thus, we have not performed a sufficiently detailed analysis. We
have not fully exploited the fact that the $x$-quadratic in
\eqref{eq:aijexp} must vanish for all $i, j$. We use the fact that
the relations in Lemma \ref{lem:sixj} must hold for all $j$, and
substitute for different choices of $i$ and $j$ in
\eqref{eq:aijexp}.

There are two cases: for all $j \neq k_0$ we have equally signed
terms, and for some $j_0 \neq k_0$ we have oppositely signed
terms. The contribution from all terms being equally signed is at
most $\frac{9n \cdot 2^{n-1}}{2^{6n}}$; this follows immediately
from there being 2 choices for the $x$-tuples for each $j \neq
k_0$.

Assume for some $j_0$ that we have oppositely signed terms; for
definiteness, say $x_{2,j_0} = -x_{3,j_0} = x_{5,j_0} =
-x_{6,j_0}$ (and of course $x_{1,j_0} = x_{4,j_0}$). From
\eqref{eq:aijexp} we have \begin{equation}x_{1,i}x_{1,j_0} +
x_{2,i}x_{2,j_0} + x_{3,i}x_{3,j_0} - x_{4,i}x_{4,j_0} -
x_{5,i}x_{5,j_0} - x_{6,i}x_{6,j_0}. \end{equation}We substitute
in the values for the $x$'s at $j_0$. Note that $x_{1,i} =
x_{4,i}$, so $x_{1,i}x_{1,j_0} - x_{4,i}x_{4,j_0} = 0$. We find
\begin{equation}\label{eq:sixfourpiece} x_{2,j_0} \cdot (x_{2,i} - x_{3,i} -
x_{5,i} + x_{6,i}) \ = \ 0; \end{equation}however, the tuple
$(x_{2,i}, x_{3,i}, x_{5,i}, x_{6,i})$ must satisfy one of the
relations in Lemma \ref{lem:sixj}.

A priori, all of the six possibilities in Lemma \ref{lem:sixj}
should be available to this tuple. If we are in the case of an
equally signed term, then \eqref{eq:sixfourpiece} is satisfied.
If, however, the tuple is oppositely signed, then one of the two
possibilities leads to a contradiction (i.e., an $x$-sum is
non-zero, and hence an $a$-sum will vanish; this would not
necessarily be the case if $m = 4$). Namely, if the second case
occurs and $x_{2,i} = -x_{3,i} = -x_{5,i} = x_{6,i}$, then the
$x$-sum in \eqref{eq:sixfourpiece} is non-zero. Thus this case
cannot occur, and for indices $i \neq k_0, j_0$, there are only $2
\cdot 4$ possibilities for the tuples, and not $2 \cdot 6$ (there
are two possibilities from $x_{1,i} = x_{4,i}$; then we saw of the
six possibilities for the rest, only four work). There are
$n(n-1)$ ways (order matters) to choose two indices $j_0, k_0$
(and for $k_0$, there are 9 ways to choose the matchings). For the
index $j_0$, there are 2 different structures of oppositely signed
terms. Each structure is determined by $x_{2,j_0}$ (two choices);
there are also two choices for $x_{1,j_0}$. Thus for $j_0$ there
is a contribution factor of 8. For the remaining $n-2$ indices,
each gives rise to 8 tuples. Each such tuple has $x_{1,1} \cdots
x_{6,n} = 1$, and the sum contributes $|\mathcal{F}|$.

Recall we divide the average by $2^{6n}$, the number of tuples.
The contribution from condition (2) holding for some index $k_0$
and at least one index $j_0$ is oppositely signed terms is
\begin{equation}\le \ 9 \cdot 8 \cdot n(n-1) \cdot
\frac{8^{n-2}}{2^{6n}} \ = \ \frac{9n(n-1)}{8} \frac{1}{2^{3n}};
\end{equation}the total contribution from condition (2) holding at least
once is therefore at most \begin{equation}\frac{9n(n-1) + 9n
2^{2-2n}}{8} \frac{1}{2^{3n}}. \end{equation} Note if condition
(3) holds for some index $k_0$, by changing each $x_{i,k_0}$ to
$-x_{i,k_0}$, then condition (2) holds. Thus the contribution from
condition (3) holding is also at most
$\frac{9n(n-1)+9n2^{2-2n}}{8} \frac{1}{2^{3n}}$. Similarly,
condition (4) holding is equivalent to condition (1) holding by a
change of variable. If condition (1) or (4) holds for each index
$i$, assuming such terms contribute fully, there are at most $2^n$
such tuples, giving a contribution bounded by
$\frac{2^n}{2^{6n}}$. Adding these bounds completes the proof of
Theorem \ref{thm:6thnonhomog}. \qed

\begin{rek} The above analysis was greatly simplified by the
presence of the linear terms in the polynomial $f(x)$. Without
relations \eqref{eq:bkexp} and \eqref{eqsixxs}, the analysis would
be significantly more involved. \end{rek}

\section{Bounds for $n \leq 10$ variables} \label{smallvarsection}

In this section, we prove upper bounds on the norm of $\tilde
S(f)=\tilde S(f,n,m)$ for  $n \leq 10$ and arbitrary odd modulus
$m \geq 3$. We shall sometimes call $\tilde S(f)$ ``the
exponential sum for polynomials of $n$ variables''.  When no
ambiguity results, we write $\tilde S$ instead of $\tilde S(f)$
(particularly for $n=3$ and $n=5$).

\begin{thm}
Let $f$, $n$, $q$, $S$ be as defined in Section \ref{sec:defnot},
and suppose $n\leq 10$. Then \begin{equation}|S| \ = \
2^{-n}|\tilde S|\ \leq \ \left({q\over
2}\right)^{\lfloor{{n+1}\over 2}\rfloor}.\end{equation}\end{thm}

\textbf{Proof.} It follows  from Lemma 3.5 of Green \cite{green}
(which easily generalizes to arbitrary odd moduli) that it is
sufficient to prove this for odd $n$ less than 10. We will first
dispose of some easy cases when the number of variables is 1 or 2,
and also when the graph $G$ has no vertex of degree 2 or more.  We
then consider in detail what happens when $n=3,5,7$, and 9.

The idea is that unless the polynomial $f$ has a special form, we
will be able to prove very small upper bounds on $|S(f)|$, which
we use in turn to prove bounds on the normalized sum for
polynomials in larger numbers of variables.

A key ingredient in the proof is the fact that $\cos(k\theta)$ is
a polynomial of degree $k$ in $\cos\theta$; these are the classic
Chebyshev polynomials.  We will use these in a slightly altered
form: $2\cos(k\theta)= Q_k(2\cos\theta)$, where the polynomials
$Q_k$ are given by the recurrence
\begin{equation}\begin{split}
Q_0(x) & \ = \  2 \\
Q_1(x)& \ = \ x \\
Q_{k+1}(x)&\ = \ xQ_k(x)-Q_{k-1}(x).\end{split} \end{equation}

We will often also need to prove that for some univariate
polynomial $g$, $g(q)>0$.  This will always follow from the fact
that $g$ is positive on the half-open interval $[\sqrt 3,2)$.
Whenever this is the case, the claim can easily be verified by
elementary calculus, but we will omit this verification in the
argument below, and simply assert $g(q)>0$.

\bigskip
\noindent{\it Case 1: $n=1$.} In this case
\begin{equation}\begin{split} f(x)& \ = \  ax \\ \tilde
S & \ = \  \omega^a-\omega^{-a},\end{split}\end{equation} so
\be|\tilde S|\ \leq\ q\end{equation}and \begin{equation}|S|\ \leq\
{q\over 2},\end{equation}as required, with equality if and only if
$a=\pm c$.

It is interesting to see what happens if $a$ is not $\pm c$. In
this case, we actually find \begin{equation}|S|\ \leq\
\left({q\over 2}\right)^9.\end{equation}To see this, we note that
$|\tilde S|$ is bounded above by $r=Q_2(q)=q^3-3q$. The claim
follows from the fact that \begin{equation}q^9-256q^3+256\cdot 3q\
\geq\ 0.\end{equation}
\bigskip
\noindent{\it Case 2: $n=2$.} While the theorem for two variables
follows from the one-variable result, we need more detailed
information for later arguments. For two variables,
\begin{equation}\begin{split} f(x,y) & \ = \  Axy+Bx+Cy \\
\tilde S & \ = \
\omega^A(\omega^{B+C}+\omega^{-(B+C)})-\omega^{-A}(\omega^{B-C}+\omega^{-(B-C)}).\end{split}
\end{equation}
If $B=C=0$ then we get the maximum value $q$ when $A=\pm c$,
giving the theorem for $n=2$.  Otherwise we find, as argued above,
$|\tilde S|\leq r<q^9/2^8$.  This gives a bound of $q^9/2^{10}$
for $|S|$.  Since $\frac12 \leq \left({q\over 2}\right)^4$, we get
a bound of $\left({q\over 2}\right)^{13}$ for $|S|$.

If either $B+C$ or $B-C$ is nonzero, then we get a bound on
$|\tilde S|$ of
\begin{equation}2+
\max_{\alpha\in\Z_m\backslash\{0\}}|\omega^{\alpha}
+\omega^{-\alpha}|\ = \ 2+s\ = \ 2+Q_2(q)\ = \
q^2.\end{equation}This bound is attained only if $A=0$ and $B=\pm
C=\pm c$, that is, with the linear polynomial $\pm cx \pm cy$. Any
other linear polynomial gives a bound of
\begin{equation}s+2\cos{{2\pi}\over m}\ = \ (q^2-2)+(q^4-4q^2+2)\
= \ q^4-3q^2\ \leq\ \left({q\over 2}\right)^{10}.\end{equation}
For a nonlinear polynomial we get a bound of
$2\omega^A+s\omega^{-A}$, attained when $B=-C=c$.  This has its
largest absolute value when $A=2c$, in which case we find
\begin{equation}\begin{split}
|\tilde S|^2&\ = \
(2\omega^{2c}+s\omega^{-2c})(2\omega^{-2c}+s\omega^{2c})\\
&\ = \ 4+2s(\omega^{4c}+\omega^{-4c})+s^2\\
&\ = \ 4+2s(s^2-2)+s^2\\
&\ = \ 4+2s^3-4s+s^2\\
& \ = \ 2q^6-11q^4+16q^2. \end{split}
\end{equation}
We verify  that for $x\in[\sqrt 3,2)$,
\begin{equation}x^5/8-\sqrt{2x^6-11x^4+16x^2}\ \geq\ 0.\end{equation} This
makes the normalized sum smaller than $\left({q\over
2}\right)^{5}$.

To summarize:  For $n=2$ we achieve the maximal value of
$\left({q\over 2}\right)$ for the magnitude of the normalized sum
when $f(x,y)=\pm cxy$.  We achieve the largest sub-maximal value
of $\left({q\over 2}\right)^{2}$ when $f(x,y)=\pm cxy\pm cxy$.  In
all other cases the magnitude of the normalized sum is less than
$\left({q\over 2}\right)^{5}$.

\bigskip
\noindent{\it Case 3. $G$ has no vertex of degree greater than 1.}
Let $n$ be any odd number of variables.  If $G$ has no vertex of
degree at least 2, then $f$ decomposes as a sum of polynomials of
degree 1 and 2 over disjoint sets of variables, and the normalized
sum $S$ for $f$ is the product of the normalized sums for each of
these polynomials.  The largest magnitude for this sum occurs when
the graph consists of $(n-1)/2$ edges and a single isolated
vertex, and when each of the associated linear and quadratic
polynomials has the largest possible normalized sum.  This
implies\begin{equation} f(x_1,\ldots,x_n)\ = \ \pm
cx_1x_2\pm\cdots \pm cx_{n-2}x_{n-1}\pm cx_n\end{equation}(up to a
permutation of the variables), giving a normalized sum whose
magnitude is  $\left({q\over 2}\right)^{{n+1}\over 2}$, as
required by the theorem.  In any other instance, the foregoing
analysis shows the normalized sum to be bounded above by
$\left({q\over 2}\right)^{{n+3}\over 2}$, which is attained when
the graph consists of three isolated vertices and $n-3$ edges.

\bigskip

\noindent{\it Case 4. $n=3$.} In this case we write
\begin{equation}\begin{split}\tilde S&\ = \ \omega^{\alpha}(\omega^{\beta}(\omega^{\gamma}-\omega^{-\gamma})-\omega^{-\beta}(\omega^{\delta}-\omega^{-\delta}))\\
& \ \ \ \ \ \ - \
\omega^{-\alpha}(\omega^{\beta'}(\omega^{\gamma'}-\omega^{-\gamma'})-\omega^{-\beta'}(\omega^{\delta'}-\omega^{-\delta'})),
\end{split}\end{equation} where
\begin{equation}\begin{split}
\alpha& \ = \ a_{12}\\
\beta&\ = \ a_1+a_2\\
\beta'&\ = \ a_1-a_2\\
\gamma& \ = \ a_{13}+a_{23}+a_3\\
\gamma'&\ = \ a_{13}-a_{23}+a_3\\
\delta & \ =  \ a_{13}+a_{23}-a_3\\
\delta' & \ = \ a_{13}-a_{23}-a_3.\end{split}
\end{equation}
We may assume with no loss of generality that $a_3\neq 0$.
 (If all the linear coefficients were zero
then $f$ would be homogeneous and $S=0$. Otherwise we can renumber
the variables to assure that $a_3$ is nonzero.)

Suppose first that all four of the subexpressions
$\omega^{\epsilon}-\omega^{-\epsilon}$ occurring in the above
equation for $S$ have the maximum possible magnitude; that is,
$\epsilon=\pm c$. If $\gamma=\delta$, we conclude (using the fact
that it is possible to divide by 2 in $\Z_m$ as $m$ is odd) that
$a_3=0$, contrary to assumption.  So $\gamma=-\delta$.  Likewise
we conclude $\gamma'=-\delta'$. This implies
$a_{13}+a_{23}=a_{13}-a_{23}=0$, so $a_{13}=a_{23}=0$.  Thus $G$
has no vertex of degree 2 or more.  By the results of the last
section we get a bound of $\left({q\over 2}\right)^{2}$ for the
normalized sum, with this largest value occurring only when $f$ is
\begin{equation}\pm cx_1x_2\pm cx_3.\end{equation}
Suppose that 3 of the 4 subexpressions in question are maximal.
This implies (up to some sign changes and renumbering of
variables): \begin{equation}a_{13}\ = \ c,\ \ \ a_{23}\ = \ -c,\ \
\ a_3\ = \ c,\end{equation}so that \begin{equation}\gamma\ =\ c,\
\ \ \delta\ = \ -c,\ \ \ \gamma'\ = \ -c,\ \ \ \delta'\ = \
-3c.\end{equation}So now \begin{equation}\tilde S\ = \
i((q\omega^{\alpha}(\omega^{\beta}+\omega^{-\beta})+
\omega^{-\alpha}(q\omega^{\beta'}+r\omega^{-\beta'})).\end{equation}If
$\alpha=\beta=\beta'=0$, then we get $|\tilde S|=3q+r=q^3$. So the
normalized sum is bounded by $q^3/8$, which is attained when $f$
has the form \begin{equation}\pm(cx_1x_3\pm cx_2x_3\pm
cx_3).\end{equation}If $\beta$ and $\beta'$ are both zero and
$\alpha$ is nonzero, we get \begin{equation}S\ = \
2qi\omega^{\alpha}+(q+r)i\omega^{-\alpha}.\end{equation}Thus
\begin{equation}\begin{split}
|\tilde S|^2&\ = \ (2q\omega^{\alpha}+(q+r)\omega^{-\alpha})(2q\omega^{-\alpha}+(q+r)\omega^{\alpha})\\
& \ = \ 4q^2+(q+r)^2+2q(q+r)(\omega^{2\alpha}+\omega^{-2\alpha})\\
&\ = \
4q^2+(q^3-2q)^2+2q(q^3-2q)(\omega^{2\alpha}+\omega^{-2\alpha}).\end{split}
\end{equation}
This is maximized when $2\alpha=1$ in $\Z_m$, which gives
\begin{equation}4q^2+(q^3-2q)^2+2q(q^3-2q)(q^4-4q^2+2).\end{equation}
We can  bound the square root of this expression on $[\sqrt 3,2)$
and find the normalized sum is less than $(q/2)^6$.  If $\beta$
and $\beta'$ are not both zero, then we get the maximal value when
$\alpha=0$ and $\beta=\beta'=2c$. The result is
\begin{equation}\tilde S\ = \
i(2q\omega^{2c}+(q+r)\omega^{-2c}),\end{equation}again giving the
bound $(q/2)^6$ for the normalized sum.

We now consider the case when no more than 2 of the subexpressions
$(\omega^{\epsilon}-\omega^{-\epsilon})$ are maximal. In this case
(remembering $a_3\neq 0$) there are no solutions for the system of
four equations in which two of the $\epsilon$  are $\pm c$ and the
other two are $\pm 3c$ (which would give a bound of $2(q+r)$).
Instead, we cannot get any value larger than
$2q+r+|\omega^{5c}-\omega^{-5c}|$.  This will happen with
$a_{13}=2c$, $a_{23}=-2c$, $a_3=c$.  We find
\begin{equation}|\omega^{5c}-\omega^{-5c}|\ = \
2\cos\left({{5\pi}\over{2m}}\right) \ = \ Q_5(q)\ = \
q^5-5q^3+5q,\end{equation}so that $|\tilde S|$ is bounded above by
\begin{equation}2q+(q^3-3q)+(q^5-5q^3+5q)\ = \ q^5-4q^3+4q.\end{equation}This implies that the normalized sum's magnitude is less than
$(q/2)^9$.

We summarize what happens in the 3-variable case. We are assuming
$a_3\neq 0$. We get the maximum magnitude for the normalized sum
of $(q/2)^2$ when $f$ is $\pm cx_1x_2\pm cx_3$.  We get the second
largest value of $(q/2)^3$ only if $f$ is either linear or has the
form $\pm(cx_1x_3\pm cx_2x_3\pm cx_3)$. In all other cases the
bound is at most $(q/2)^4$.

For future reference, it is worth thinking explicitly about the
case where $a_3=0$ and $a_{13}$, $a_{23}$ are both nonzero.  We
get $\gamma=\delta$ and $\gamma'=\delta'$.  Furthermore, we cannot
have $\gamma=\pm\delta$ without making one of $a_{12}$ or $a_{13}$
zero.  The largest norm possible occurs when $\gamma=c$ and
$\gamma'=3c$, in which case \begin{equation}\tilde S\ = \
\omega^{\alpha}(\omega^{\beta}-\omega^{-\beta})qi+\omega^{-\alpha}(\omega^{\beta'}-\omega^{-\beta'})ri,\end{equation}so
\begin{equation}|\tilde S|\ \leq\ q^2+qr\ = \ q^2+q^4-3q^2\ = \
q^4-2q^2,\end{equation}which gives a normalized sum whose
magnitude is no more than $(q/2)^6$.

\bigskip

\noindent{\it The ``General Case''.} ``General'' here means 5, 7,
or 9. Note again that if $G$ has no vertex of degree two or higher
then by Case 3 we have all the information we need (in particular,
we obtain the stated bound on the normalized sum, valid for
arbitrary $n$). Accordingly, suppose $G$ has a vertex of degree 2
or more. We may assume without loss of generality that this is
vertex $n$, and that $a_{n-1,n}$ and $a_{n-2,n}$ are both nonzero.

We write $f^{++}$, $f^{-+}$, etc. for the four $(n-2)$-variable
polynomials  formed by setting $x_1$ and $x_2$ to $\pm 1$ and then
setting the constant term of the resulting polynomial to zero. For
example, if
\begin{equation}f(x_1,x_2,x_3,x_4,x_5)=a_{12}x_1x_2+a_{23}x_2x_3+a_{34}x_3x_4+a_1x_1+a_3x_3+a_4x_4+a_5x_5,\end{equation}then \begin{equation}f^{-+}(x_3,x_4,x_5)\ = \
a_{34}x_3x_4+(a_3+a_{23})x_3+a_4x_4+a_5x_5.\end{equation}We denote
by $S^{++}$, $S^{-+}$, etc., the unnormalized sums of the
$f^{\pm\pm}$, and by $G^{\pm\pm}$  the graph (it's the same for
all four polynomials) of the $f^{\pm\pm}$. We now have
\begin{equation}S=\omega^{a_{12}}(\omega^{a_1+a_2}S^{++}+\omega^{-(a_1+a_2)}S^{--})-
\omega^{-a_{12}}(\omega^{a_1-a_2}S^{+-}-\omega^{-(a_1-a_2)}S^{-+}).\end{equation}Note
that each of the $f^{\pm\pm}$ has a vertex of degree at least 2 in
the associated graph.

We want to show that the largest possible normalized sum for
polynomials in $x_3,\ldots, x_n$ with $a_{n-1,n}$ and $a_{n-2,n}$
both nonzero, occurs {\it only}  when the polynomial has the form
\begin{equation}\pm cx_3x_4\pm cx_5x_6\pm\cdots\pm cx_{n-1,n}x_{n-2,n}\pm
cx_n\end{equation}(up to a permutation of $\{3,4,\ldots,n-3\}$).
In this case the magnitude  of the unnormalized sum for $n-2$
variables is $2^{(n-5)/2}q^{(n+1)/2}$. This would imply that the
normalized sum for polynomials in $n$ variables is bounded above
by
\begin{equation}2^{-n}\cdot 4\cdot
2^{(n-5)/2}q^{(n+1)/2}\ = \ \left({q\over 2}\right)^{{n+1}\over
2},\end{equation}as required by the theorem.  Observe that in our
study of three-variable polynomials we have already established
this claim in the case $n=5$. We proceed to show it for $n=7$ and
$n=9$. We really want to show by induction that this claim holds
for all odd $n$. Let us suppose then that this property of
polynomials in $n-2$ variables holds, and see how close we can
come to completing the inductive proof.

How many of the  $S^{\pm\pm}$ can give the optimal magnitude of
$2^{(n-5)/2}q^{(n+1)/2}$ for polynomials in $n-2$ variables with a
vertex of degree 2? Suppose first that all four of these sums are
optimal. Then by induction each of the $f^{\pm\pm}$ is
\begin{equation}\pm cx_3x_4\pm cx_5x_6\pm\cdots\pm
cx_{n-1,n}x_{n-2,n}\pm cx_n\end{equation}We thus have for $3\leq
i<n$,
\begin{equation}a_i\pm a_{1i}\pm a_{2i}\ = \ 0,\end{equation}which implies
\begin{equation}a_{1i}\ = \ a_{2i}\ = \ 0.\end{equation}We also have
\begin{equation}a_n\pm a_{1n}\pm a_{2n}\ = \ \pm c.\end{equation}If three of
the four values $a_n\pm a_{1n}\pm a_{2n}=\pm c$  are equal, we
find $a_{1n}=a_{2n}=0$ (so that all four of the values are equal),
and thus $G$ is disconnected, with $\{1,2\}$ as a separate
component. In this case $|S|$ cannot exceed the product of the
magnitudes of the  sums associated with the components, namely
\begin{equation}2^{(n-5)/2}q^{(n+1)/2}\cdot 2q\ = \
2^{(n-3)/2}q^{(n+3)/2}.\end{equation}Observe that this arises
precisely when $f$ has the form \begin{equation}\pm cx_1x_2\pm
cx_3x_4\pm cx_5x_6\pm\cdots\pm cx_{n-1,n}x_{n-2,n}\pm
cx_n.\end{equation}This gives a bound on the normalized sum of
$(q/2)^{{n+3}\over 2}$. To complete the induction we will have to
show that every other possible form for $f$ gives a strictly
smaller value.

We may thus suppose that two of the four values
\begin{equation}a_n\pm a_{1n}\pm a_{2n}=\pm c\end{equation} are $c$ and two
are $-c$. We can assume without loss of generality that
\begin{equation}a_n+a_{1n}+a_{2n}\ = \ c.\end{equation}If we also have \begin{equation}a_n-a_{1n}-a_{2n}\ = \ c,\end{equation}then $a_n=c$ and $a_{1n}+a_{2n}=0$.
This would imply that both $\pm(a_{1n}-a_{2n})$ equal $-c$, which
is impossible. Thus \begin{equation}a_n-a_{1n}-a_{2n}\ = \
-c,\end{equation}which implies $a_n=0$ and $a_{1n}+a_{2n}=c$. This
implies $a_{13}-a_{23}=\pm c$, and thus either $a_{13}=0$ or
$a_{23}=0$. The result is that \begin{equation}\tilde S\ = \
S^{++}\left[\omega^{a_{12}}(\omega^{a_1+a_2}
-\omega^{-(a_1+a_2)})\pm\omega^{-a_{12}}
(\omega^{a_1-a_2}-\omega^{-(a_1-a_2)})\right].\end{equation}The
largest possible magnitude for the bracketed expression is $q^2$,
giving a bound of $q^2\cdot 2^{(n-5)/2}q^{(n+1)/2}$ for $|\tilde
S|$, and thus of $(q/2)^{{n+5}\over 2}$ for $|S|$.

We now suppose that exactly three of the $S^{\pm\pm}$ have
magnitude $2^{(n-5)/2}q^{(n+1)/2}$. Note that whenever at least
one of the $S^{\pm\pm}$ has this form, the graph $G^{\pm\pm}$ is
disconnected, with a component consisting of the vertices
$\{n-2,n-1,n\}$. Thus each $S^{\pm\pm}$ is the product of the sum
$S_3^{\pm\pm}$ associated with some three-variable polynomial
$f_3^{\pm\pm}$ and the sum associated with an $(n-5)$-variable
polynomial.  By the inductive hypothesis, the sum for an
$(n-5)$-variable polynomial has magnitude bounded above by
$(q/2)^{{(n-5)}\over 2}$.

We can suppose without loss of generality that the three optimal
sums  are $S^{++}$, $S^{+-}$, and $S^{-+}$. We again find
\begin{equation}a_{1,n-1}\ = \ a_{2,n-1}\ = \ a_{1,n-2}\ = \
a_{2,n-2}\ = \ 0.\end{equation}We also have
\begin{equation}\begin{split} a_n+a_{1n}+a_{2n}& \ = \
\pm c, \\ a_n+a_{1n}-a_{2n}& \ = \ \pm c,\\
a_n-a_{1n}+a_{2n}& \ = \ \pm c.\end{split}\end{equation} If all
three right-hand sides above are equal, we again get
$a_{1n}=a_{2n}=0$, which will put us back in the previous case. If
the first two right-hand sides are equal, and the third is
opposite, we find $a_n=a_{2n}=0$, which again puts us back in the
previous case.  We may thus suppose  that the first right-hand
side is $c$, so that the second is $-c$. We then obtain
\begin{equation}a_n\ = \ -c,\ \ \ a_{1n}\ = \ a_{2n} \ = \
c.\end{equation} Thus $|S_3^{++}|=|S_3^{+-}|=|S_3^{-+}|=q^3$, and,
as we found in the section on 3 variables, $|S_3^{--}|$ is the
magnitude of the sum for the 3-variable polynomial
$cx_1x_3-cx_2x_3-3cx_3$. We find, reasoning as in the section on
three variables, that this is
\begin{equation}q+2r+q^5-5q^3+5q\ = \ q^5-3q^3.\end{equation} Thus the sum
of the $|S_3^{\pm\pm}|$ is no more than $q^5$, so that
\begin{equation}2^{-n}|S|\ \leq\ (q/2)^5\cdot (q/2)^{{(n-5)}\over
2}\ = \ (q/2)^{{n+5}\over 2}.\end{equation} In the case where one
or two of the $S_3^{\pm\pm}$ have the value $q^3$, the same
reasoning applies and leads to a bound (not the best possible!) of
$(q/2)^{{n+5}\over 2}$ for the normalized sum.

We are thus left with the case where none of the $S_3^{\pm\pm}$
attain the maximal value $q^3$. In this instance we can no longer
suppose that $\{n-2,n-1,n\}$ forms a separate component of
$G^{\pm\pm}$, so we will have to be content to argue for specific
values of $n$.

For $n=5$, the analysis of the the 3-variable case shows that each
$|S^{\pm\pm}|$ is bounded above by $q^6/8$, which by the triangle
inequality gives the bound $q^6/2$ for $|\tilde S|$. This, in
combination with the calculations above, shows that if $f$ is a
polynomial in 5 variables such that $G$ has a vertex of degree at
least 2, and $f$ is not of the special form \begin{equation}\pm
cx_1x_2+\pm cx_3x_4\pm cx_5x_6\pm\cdots\pm cx_{n-1,n}x_{n-2,n}\pm
cx_n,\end{equation}then $|S|\leq q^5$.  This allows us to extend
our ``induction'' to seven variables:  If $f$ is a polynomial on 7
variables for which $G$ has a vertex of degree at least 2, either
$G$ has the special form above, or $|S|$ is bounded above by
$4q^5$.  Applying the argument one more time shows that for
polynomials in 9 variables, in all cases we get a bound on $|S|$
of $16q^5$, which gives a bound on the $|S|$ of $(q/2)^5$, as
required.

\begin{rek} Where do things fall apart?  Observe that the induction
fails precisely when none of the $S_3^{\pm\pm}$ are maximal (for
polynomials whose graphs have a vertex of degree at least 2). We
made use of the fact that if one of the $S_3^{\pm\pm}$ is maximal
in this sense, then $G^{\pm\pm}$ has a component with three
vertices, and this condition is sufficient for the induction to
carry through. Ironically, the principal obstruction to completing
the proof occurs for polynomials whose sums we expect to have
values that are very far from the conjectured upper
bound.\end{rek}

\section{Fourier Bounds} \label{fouriersection}

In this section, we use Fourier analytic methods to provide bounds
for  $S(f)$, where $f$ is a polynomial in $Z_m^2[n]$ whose graph
$G(f)$ is (almost) acyclic (the precise definition is given
below). We first need to establish some notation.

\subsection{Notation}
Let $\Omega=\{1, -1\}$ and define $L^2=L^2(\Omega^n)= \{g\ |\ g:
\Omega^n \rightarrow \mathbb{C}\}$.  Let $[n]$ denote the set
$\{1,2, \ldots, n\}$. The set of functions $\chi_S \in L^2$ for $S
\subseteq [n]$ where \begin{equation}\chi_S(x_1, x_2, \ldots,
x_n)\ =\ \prod_{i \in S} x_i \end{equation}form an orthogonal
Fourier basis for $L^2$ where the inner product  of  functions $f$
and $g$ is defined as follows: \begin{equation}\langle f, g
\rangle\ =\ \sum_{y \in \Omega^n}  f(y) \overline{g(y)}.
\end{equation}where $\overline{z}$ is the complex conjugate of $z
\in \mathbb{C}$.

Thus any function $g \in L^2$ can be written as \begin{equation}g\
=\ \sum_{S \subset [n]} c_S(g) \chi_S \end{equation}which we call
the Fourier expansion of $g$, where $c_S(g)$ is a particular
Fourier coefficient in the expansion.

Since the $\{\chi_S | S \subset [n]\}$ is an orthogonal basis, we
can express $c_S$ as follows: \begin{equation}c_S = \langle g,
\chi_S \rangle = \sum_{y \in \{1, -1\}^n} g(y)
\overline{\chi_S(y)} = \sum_{y_i \in \{1,-1\}^n} \left(\prod_{i
\in S}  y_i \right) \ g(y_1, y_2, \ldots, y_n). \end{equation}This
implies that the exponential sum $S(f)$ under consideration is the
Fourier coefficient $c_S(g)$ when $S=\{1,2,\ldots, n\}$ and
$g=\omega^{f(x_1, x_2, \ldots, x_n)} \in L^2$. We let
$\hat{c}_S(f)= c_S(\omega^{f})$, which we sometimes denote as
$\hat{c}_S$ when $f$ is obvious from the context. Our goal then is
to prove that $\hat{c}_{[n]}(f)$ is exponentially small for every
polynomial $f \in \Z_m^2$.

It is possible, in some cases, to give an explicit computation of
the Fourier expansion, which we now show.  Let $f(x_1, x_2,
\ldots, x_n)
 = \sum_{i \not= j} a_{ij} x_i x_j + \sum_i a_i x_i$ be a quadratic
polynomial of $n$ variables where $a_{ij}, a_i \in \mathbb{Z}_m$.
Observe that \begin{equation}\omega^{a_{ij} x_i x_j}\ = \
\frac{1}{2}\, (\om^{a_{ij}} - \om^{-a_{ij}})\, x_i x_j +
\frac{1}{2}\, (\om^{a_{ij}} +\om^{-a_{ij}}) \end{equation}and
\begin{equation}\omega^{a_i x_i} \ = \ \frac{1}{2}\, (\om^{a_{i}}
- \om^{-a_{i}})\, x_i  + \frac{1}{2}\, (\om^{a_{i}} +\om^{-a_{i}})
\end{equation}since $x_i, x_j \in \{1,-1\}$. We set $\lambda(x)= (\omega^x -
\omega^{-x})/2$ and $\mu(x) = (\omega^x + \omega^{-x})/2$. Thus we
are interested in  the coefficient of $x_1 x_2 \ldots x_n$ when we
expand  and simplify \begin{equation}\prod_{i \not=j}
(\lambda(a_{ij})\, x_i x_j + \mu(a_{ij}))\, \prod_{i}
(\lambda(a_i) x_i + \mu(a_i)), \end{equation}using the relations
$x_i^2=1$ for all $1 \leq i \leq n$.

\subsection{Bounds  on Fourier Coefficients for a special class of polynomials}

 Recall that for a polynomial $f(x_1, x_2, \ldots, x_n)$ we can associate the weighted undirected graph $G=G(f)=(V, E)$
with vertices $V=\{1, 2, \ldots, n \}$ and edge set $E= \{\{i,
j\}|\, a_{ij} \not= 0\}$, where edge $\{i, j\}$ has weight
$a_{ij}$ (when $a_{ij} \not= 0$). We now show that when $G(f)$ is
a tree, \emph{every} Fourier coefficient is small.

\begin{lem} \label{allbound}
If $G(f)$ is a tree with $n$ vertices where $n \geq 2$, then
$|\hat{c}_{S}(f)| \leq
\left(\cos\left(\frac{\pi}{2m}\right)\right)^{n-1}$ for all $S
\subseteq [n]$.
\end{lem}

\textbf{Proof.} The bound holds when $n=2$ (see proof of Theorem
\ref{grandmother} (ii)).

Now let $f(x_1, \ldots, x_n)$ be such that $G(f)$ is a tree with
$n$ vertices where $n > 2$. Let $\{i,j\}$ be an edge in $G(f)$
with weight $a_{ij}$ such that $j$ is a leaf.  Set $f = f' +
a_{ij} x_i x_j + a_j x_j$ where $f'$ is independent of $x_j$.

Since \begin{equation}\omega^{f}\ =\ \omega^{f'} \left(
\h{\lambda(a_{ij})} x_i x_j + \h{\mu(a_{ij})} \right) \left(
\h{\lambda(a_j)} x_j + \h{\mu(a_j)} \right), \end{equation}the
coefficient $\hat{c}_S(f)$ can be written in terms  of the Fourier
coefficients $\hat{c}(f')$. Then for any $S \subseteq ([n]
\setminus \{j\})$,
 \begin{equation}\hat{c}_S(f)\ =\  \frac{\mu(a_{ij})}{2} \frac{\mu(a_j)}{2} \hat{c}_{S}(f') + \frac{\lambda(a_{ij})}{2} \frac{\lambda(a_j)}{2}\, \hat{c}_{S \symd \{i\}}(f')  \end{equation}where $\symd$  refers to the symmetric difference of two sets: $
\symd B= (A \setminus B) \cup (B \setminus A)$. Similarly for any
subset $S \subseteq [n]$ such that $j \in S$,
\begin{equation}\hat{c}_S(f)\ = \  \frac{\mu(a_{ij})}{2}
\frac{\lambda(a_j)}{2}\, \hat{c}_{S \symd \{j\}}(f') +
\frac{\lambda(a_{ij})}{2} \frac{\mu(a_j)}{2} \hat{c}_{S \symd
\{i\}}(f') \end{equation} Assume (via induction on $n$) that
$|\hat{c}_S(f')|\ \leq\ (\cos\left(\frac{\pi}{2m})\right)^{n-2}$.
Then
\begin{equation}|\hat{c}_S(f)|\ \leq\  \frac{1} {4}
\left(\cos\left(\frac{\pi}{2m}\right) \right)^{n-2} \left( \left|
\mu(a_{ij}) \mu(a_j) \right| + | \lambda(a_{ij}) \lambda(a_j)|
\right) \end{equation}when $j \in S$ and
\begin{equation}|\hat{c}_S(f)|\ \leq\ \frac{1} {4}
\left(\cos\left(\frac{\pi}{2m}\right)\right)^{n-2} \left(
|\mu(a_{ij})\lambda(a_j)| + |\lambda(a_{ij}) \mu(a_j)| \right)
\end{equation}when $j \not\in S$.

We first consider the case when $j \not\in S$ (the other case is
handled similarly). If $a_{ij}=a_j$, then
\begin{equation}\begin{split} |\hat{c}_S(f)|&\ \leq\  \frac14
(|\hat{c}_{S \symd \{i\}}(f')|+ |\hat{c}_{S \symd \{j\}}(f')|)\\ &
\ \leq \ \frac12 \left(\cos\left(\frac{\pi}{2m}\right)
\right)^{n-2} \ \leq\ \left(\cos \left(\frac{\pi}{2m}\right)
\right) ^{n-1}.\end{split}\end{equation} If $a_{ij} \not= a_j$,
\begin{equation}\left| \mu(a_{ij}) \mu(a_j) \right| + |
\lambda(a_{ij}) \lambda(a_j)|  = 4 (|\sin(\theta)| |\sin(\alpha)|
+ |\cos(\theta)| |\cos(\alpha)|) \end{equation}where $\theta= 2
\pi a_{ij}/m$ and $\alpha= 2 \pi a_j/m$ are both multiples of $2
\pi/m$.

Observe that we may reflect $\om^{a_{ij}}$ and $\om^{a_j}$ to the
first quadrant since this operation does not change the absolute
value of either the sine or cosine of their arguments. After this
transformation, $\theta$ and $\alpha$ are integral multiples of
$\pi/2m$ and are both $< \pi/2$. This implies that
\begin{equation}|\sin(\theta)| |\sin(\alpha)| + |\cos(\theta)|
|\cos(\alpha)| \ = \ \cos(\theta-\alpha).\end{equation} Since
$\theta-\alpha$ is an integral multiple of $\pi/2m$ and $\theta
\not= \alpha$ (since $a_{ij} \not= a_j$)
\begin{equation}|\hat{c}_S(f)|\ \leq\
\left|\cos\left(\frac{\pi}{2m}\right) \right|^{n-2}
|\cos\left(\frac{\pi a}{2m}\right)| \end{equation}for some $a
\not=0$, when $j \not\in S$, from which we can conclude that
$|\hat{c}_S(f)| \leq \left(\cos(\pi/2m)\right)^{n-1}$ since
$|\cos(a \pi/2m)| \leq \cos(\pi/2m)$ for all $a \not=0$.
Similarly, when $j \in S$,
\begin{equation}|\hat{c}_S(f)|\ \leq\
\left(\cos\left(\frac{\pi}{2m}\right)\right)^{n-1}.
\end{equation}\qed

\begin{rek} Observe that Lemma \ref{allbound} implies our desired bound
on the exponential sum: If $G(f)$ is a tree with $n$ vertices,
$|\hat{c}_{[n]}| \leq (\cos(\frac{\pi}{2m}))^{n-1}$.  If $G(f)$ is
a forest of disjoint trees $T_1 \cup T_2 \ldots \cup T_k$, then
$\hat{c}_S = \prod_{i=1}^k \hat{c}_{S_i}$ where $S_i$ is
restricted to vertices in $T_i$ and $S=\cup S_i$. This implies the
bound holds for a forest of trees. \end{rek}

\textbf{Proof of Theorem \ref{grandmother} (iii).} Suppose $G(f)$
is a tree and we now add a term $a_{ij} x_i x_j$ to $f$
(equivalently, add an edge of weight $a_{ij}$ to $G(f)$ between
$i$ and $j$), where we assume that there was no such term in $f$
before (if there was, this operation just modifies the weight).
Set $f'=f + a_{ij} x_i x_j$. Then, for any $S \subseteq [n]$,
\begin{equation}\hat{c}_S(f')\ =\ \lambda(a_{ij}) \hat{C}_{S \symd
\{x_i,x_j\}} + \mu(a_{ij}) \hat{c}_S(f) .\end{equation}This
implies that \begin{equation}\begin{split} |\hat{c}_S(f')|& \ \leq
\  (\max_S |\hat{c}_S(f)|)\ (|\lambda(a_{ij}) + |\mu(a_{ij})|)
\\ & \ =\  \max_S |\hat{c}_S(f)|\ (|\sin(\theta)| +
|\cos(\theta)|),\end{split}
\end{equation} where $\theta=2 \pi a_{ij}/m$ (where $\theta \not=
0, \pi/2$). Since the maximum value of $|\sin(\theta)| +
|\cos(\theta)|$ is $\sqrt{2}$, we have \begin{equation}\max_S
\size{\hat{c}_S(f')}\ \leq\ \sqrt{2} \
\max_{S}\size{\hat{c}_S(f)}.\end{equation} Clearly the same bound
holds if we add a linear term $a_i x_i$ that did not exist before.
So if $k$ such new edges are added to $G(f)$,
\begin{equation}\size{\hat{c}_S(f')}\ \leq\ 2^{k/2} \max_S
\size{\hat{c}_S(f)}\ \leq\ 2^{k/2}
\left(\cos\left(\frac{\pi}{2m}\right) \right)^{n-1}
\end{equation}Therefore when $k \leq (n-2)
\log(\frac{1}{\cos(\pi/(2m))})$, we have
 \begin{equation}\size{\hat{c}_S(f')}\ \leq\ (\cos(\pi/2m))^{n/2} \end{equation}thus obtaining the conjectured bound.

Thus if there exist a set of at most  $(n-2) \log(1/\cos(\pi/2m))$
edges from $G(f)$ whose deletion makes $G(f)$ a forest of trees,
then  \begin{equation}\size{\hat{c}_ S(f)}\ \leq\
(\cos(\pi/2m))^{n/2}
\end{equation}(recall that $q= 2 \cos(\pi/(2m))$ in the statement
of Theorem \ref{grandmother} (iii)).  \qed

\begin{rek} It is worth noting two important limitations of the above
proof:
\begin{enumerate}
\item The proof relies on a global bound for all Fourier
coefficients, whereas the only coefficient of interest is
$\hat{c}_{\{1,\dots,n\}}(f)$. \item The norm of a particular
Fourier coefficient might increase or decrease as we add
additional edges. Since we do not have the means to analyze the
behavior, we have assumed that the coefficients may increase in
norm by a factor of $\sqrt{2}$ (it is unlikely that this blowup
will occur on every edge addition and for every coefficient). A
closer analysis of this aspect might lead to a better estimate on
the number of additional edges allowed.
\end{enumerate}
\end{rek}

\section{Recent Progress and Future
Work}\label{sec:addapproachfuturework}

We believe that Conjecture \ref{holygrail} provides a tight bound
that is exponentially decreasing; while we have verified this for
$n\le 10$ and quadratic $f$, the general case is still open.

It is possible that there is more to say about sub-maximal values
of $|S(f,n,m)|$. Implicit in many of the arguments in Section
\ref{smallvarsection} is a bound on the second largest value of
$|S(f,n,m)|$. In particular, we make the following (stronger)
conjecture:

\begin{conj}[Stronger form of Conjecture \ref{holygrail}]
Let $m \geq 3$ be odd and let $n$ be a non-negative integer. Then
for quadratic $f$, \begin{equation}\size{S(f,n,m)}\ \leq\ \left( q
\over 2 \right)^{\lfloor \frac{n+1}{2} \rfloor },
\end{equation}and moreover, if $\size{S(f,n,m)} < \left( q \over 2
\right)^{\lfloor \frac{n+1}{2} \rfloor }$, then
\begin{equation}\size{S(f,n,m)}\ \leq\   \left( q \over 2
\right)^{\lfloor \frac{n+1}{2}\rfloor + 1 }.
\end{equation}\end{conj}

\begin{rek} This stronger form  has also been verified for $m=3$ by
\cite{green2} and born out by experimental evidence for small $n,
m$. \end{rek}

Lastly, we note that the  problem of bounding $S(f,n,m)$ for
polynomials $f$ of degree $2$ is only a first step. The goal is to
prove exponentially small upper bounds for all $f$ of degree
$O((\log n)^c)$ where $n$ is the number of variables. The moment
analysis can readily be carried out for such polynomials. We again
obtain square-root cancellation on average when $n \ge
\text{deg}(f)+1$, and if $\gamma<1$ is quite close to $1$ then all
but an exponentially small (in $n$) proportion of the $|S(f,n,m)|$
are bounded by $\gamma^n$.

Since the submission of this paper the fundamental problem of
proving an exponentially decreasing upper bound for $|S(f,n,m)|$
with $f$ a polynomial of fixed degree $d$ and any $n$ and $m$ has
been solved by Bourgain \cite{bourgain}, though the bounds
obtained are larger than what we feel is the true story (and for
quadratic $f$ with $m$ odd and $n\le 10$, larger than the bounds
which we show are sharp).

\section*{Acknowledgements}
We thank Avner Ash, David M. Barrington, Ron Evans, Frederic
Green, Rob Gross, John Hsia, Gene Luks and Eitan Sayag for many
enlightening conversations, and Jean Bourgain for sharing his
preprint.


\appendix

\section{Bounds when $m=3$ and $d=2$}\label{sec:appm3d2}

When $m=3$ and $d=2$, we may write \eqref{eq:sfnm} (see also
\eqref{eq:Sgnm}) as \begin{equation}\label{eq:sfnm=3} S(f,n,3) \ =
\ \frac{1}{2^n} \sum_{x_1 = -1}^{1} \cdots \sum_{x_n = -1}^{1}
\jst{x_1\cdots x_n} e_3\left( g(x) \right). \end{equation}The
presence of the Legendre symbol, coming from the factor $x_1\cdots
x_n$, complicates the arguments, giving us a mixed (additive and
multiplicative characters) complete exponential sum. We can remove
the Legendre factor by using the following identity: for $y \in
\{-1,0,1\}$, \begin{equation}\jst{y} \ = \ \frac{e_3(y) -
e_3(-y)}{i\sqrt{3}} \ = \
\begin{cases}\ \ 1 & \text{if $y = 1$}\\ \ \ 0 &\text{if $y = 0$} \\ -1 & \text{if
$y = -1$;} \end{cases} \end{equation}thus we may replace the
Legendre symbol with a product of exponentials. While this
identity can be used for any modulus (and we could use it directly
on $x_1\cdots x_n$ without passing through Legendre symbols), it
is useful only when $m=3$.

It would be natural to replace $\jst{x_1\cdots x_n}$ with
$\frac{e_3(x_1\cdots x_n) - e_3(-x_1\cdots x_n)}{i\sqrt{3}}$;
unfortunately, this would replace $S(f,n,3)$ with two exponential
sums $S'(f_1,n,3)$ and $S'(f_2,n,3)$, with $f_i$ of degree $n$
(note these sums are not mixed, composed solely of additive
characters). As Deligne's and others' bounds are of the form
$(\deg f_i - 1)^n 3^{n/2}$, this increases the degree too much to
be useful. A better approach is to let $\sigma$ be any permutation
of $\{1,\dots,n\}$ (for simplicity we consider $n$ even) and to
write \begin{equation}\jst{x_1 \cdots x_n} \ = \ \prod_{j=1}^{n/2}
\frac{e(x_{\sigma(2j)} x_{\sigma(2j-1)}) - e(-x_{\sigma(2j)}
x_{\sigma(2j-1)})}{i\sqrt{3}}. \end{equation}Expanding the product
gives $2^{n/2}$ degree $2$ exponential terms, as well as a factor
of $\left(\frac{1}{i\sqrt{3}}\right)^{n/2}$. Substituting this
into \eqref{eq:sfnm=3} yields $2^{n/2}$ complete exponential sums
$S'(f_{i,\sigma},n,3)$, where each $f_{i,\sigma}$ is of degree
$2$. If for each $f_i$ we have the homogeneous part of highest
degree is non-singular modulo $3$, then by Deligne's bound
$|S'(f_{i,\sigma},n,3)| \le \frac{3^{n/2}}{2^n}$ (recall we are
dividing by $2^n$ and not $3^n$, as initially each $x_i \in
\{-1,1\}$). Therefore for $n$ even, \begin{equation}|S(f,n,3)| \
\le \ \frac{1}{\sqrt{3}^{n/2}} \sum_{j=1}^{2^{n/2}}
|S'(f_{i,\sigma},n,3)| \ \le \ \frac{2^{n/2}}{\sqrt{3}^{n/2}}
\cdot \frac{3^{n/2}}{2^n} \ = \ \left(
\frac{\sqrt{3}}2\right)^{n/2}. \end{equation} We have shown

\begin{thm}\label{thm:sfnm=3sigma}
Let $f$ be a quadratic polynomial such that there is some
permutation $\sigma$ of $\{1,\dots,n\}$ for which the homogeneous
part of highest degree of each $f_{i,\sigma}$ is non-singular
modulo $3$. Then if $n$ is even, Conjecture \ref{holygrail} is
true for this $f$ and $m=3$.
\end{thm}

To handle odd $n$, as we must keep all the factors of degree $2$
the last factor is
$\frac{e_3(x_{\sigma(n)})-e_3(-x_{\sigma(n)})}{i\sqrt{3}}$. A
similar argument yields Conjecture \ref{holygrail} for odd $n$,
but with a slightly weaker bound, namely $\left(
\frac{\sqrt{3}}2\right)^{\lfloor n/2\rfloor}$.

To complete the investigation of $m=3$ and $d=2$ we must analyze
which $f$ satisfy the conditions of Theorem \ref{thm:sfnm=3sigma}.
For $n$ even, there are $(n-1)!!$ choices for $\sigma$ which lead
to different exponential products (the number of ways to pair $n$
objects where order does not matter); all we need is one valid
choice. As the conjecture is already known in this case, we
content ourselves with the above observation.


\bibliographystyle{plain}






\end{document}